\documentclass[smallextended,envcountsect]{svjour3}%
\usepackage{graphicx}
\usepackage{amssymb}
\usepackage{amsmath}
\journalname{}

\newcommand{\R}{{\mathbb R}}

\newcommand{\cl}{\operatorname{cl}}

\newcommand{\aff}{\operatorname{aff}}

\newcommand{\fix}{\operatorname{fix}}
\newcommand{\gph}{\operatorname{gph}}
\newcommand{\inte}{\operatorname{int}}

\newcommand{\ri}{\operatorname{ri}}

\begin{document}

\title{A Ky Fan minimax inequality for quasiequilibria on finite dimensional spaces}


\author{Marco Castellani   \and  Massimiliano Giuli \and Massimo Pappalardo}

\institute{Marco Castellani \and Massimiliano Giuli \at
             Department of Information Engineering, Computer Science and Mathematics\\
             Via Vetoio, 67100 Coppito (AQ), Italy
           \and
             Massimo Pappalardo, Corresponding author  \at
              Department of Computer Science\\
              Largo B.Pontecorvo 3, 56127 Pisa, Italy\\
              massimo.pappalardo@unipi.it
}


\maketitle

\begin{abstract}
Several results concerning existence of solutions of a quasiequilibrium problem defined on a finite dimensional space are established.
The proof of the first result is based on a Michael selection theorem for lower semicontinuous set-valued maps which holds in finite dimensional spaces.
Furthermore this result allows one to locate the position of a solution.
Sufficient conditions, which are easier to verify, may be obtained by imposing restrictions either on the domain or on the bifunction.
These facts make it possible to yield various existence results which reduce to the well known Ky Fan minimax inequality when the constraint map is constant and the quasiequilibrium problem coincides with an equilibrium problem.
Lastly, a comparison with other results from the literature is discussed.
\end{abstract}


\section{Introduction}

In \cite{Fa72} the author established the famous Ky Fan minimax inequality which concerns the existence of solutions for an inequality of minimax type that nowadays is called in literature ``equilibrium problem''.
Such a model has gained a lot interest in the last decades because it has been used in different contexts as economics, engineering, physics, chemistry and so on (see \cite{BiCaPaPa13} for a recent survey).

In these equilibrium problems the constraint set is fixed and hence the model can not be used in many cases where the constraints depend on the current analyzed point.
This more general setting was studied for the first time in the context of impulse control problem \cite{BeGoLi73} and it has been subsequently used by several authors for describing a lot of problems that arise in different fields: equilibrium problem in mechanics, Nash equilibrium problems, equilibria in economics, network equilibrium problems and so on.
This general format, commonly called ``quasiequilibrium problem'', received an increasing interest in the last years because many theoretical results developed for one of the abovementioned models can be often extended to the others through the unifying language provided by this common format.

Unlike the equilibrium problems which have an extensive literature on results concerning existence of solutions, the study of quasiequilibrium problems to date is at the beginning even if the first seminal work in this area was in the seventies \cite{Mo76}.
After that, the problem concerning existence of solutions has been developed in some papers \cite{AlRa16,Au93,AuCoIu17,CaGi15,CaGi16,Cu95,Cu97}.
Most of the results require either monotonicity assumptions on the equilibrium bifunction or upper semicontinuity of the set-valued map which describes the constraint.
Whereas other authors provided existence of solutions avoiding any monotonicity assumption and assuming lower semicontinuity of the constraint map and closedness of the set of its fixed points.

Aim of this paper is to establish several results concerning existence of solutions of a quasiequilibrium problem defined on a finite dimensional space which comes down to the Ky Fan minimax inequality in the classical setting.
Our approach is based on a Michael selection result \cite{Mi56} for lower semicontinuous set-valued maps.
Moreover the proof of our results allow one to locate the position of a solution.
The paper is organized as follows.
Section 2 is devoted to recall the results about set-valued maps which are used later.
In Section 3 we prove the main theorem and we furnish more tractable conditions on the equilibrium bifunction which guarantee that our result holds true.

\section{Basic concepts}

Let $\Phi:X\rightrightarrows Y$ be a set-valued map with $X$ and $Y$ two topological spaces.
The graph of $\Phi$ is the set
$$\gph\Phi:=\{(x,y)\in X\times Y:y\in \Phi(x)\}$$
and the lower section of $\Phi$ at $y\in Y$ is
$$\Phi^{-1}(y):=\{x\in X:y\in \Phi(x)\}.$$
The map $\Phi$ is said to be lower semicontinuous at $x$ if for each open set $\Omega$ such that $\Phi(x)\cap\Omega\ne\emptyset$ there exists a neighborhood $U_x$ of $x$ such that
$$\Phi(x')\cap\Omega\ne\emptyset,\qquad\forall x'\in U_x.$$
Notice that a set-valued map with open graph has open lower sections and, in turn, if it has open lower sections then it is lower semicontinuous.

A fixed point of a function $\varphi:X\rightarrow X$ is a point $x\in X$ satisfying $\varphi(x)=x$.
A fixed point of a set-valued map $\Phi:X\rightrightarrows X$ is a point $x\in X$ satisfying $x\in\Phi(x)$.
The set of the fixed point of $\Phi$ is denoted by $\fix\Phi$.
One of the most famous fixed point theorems for continuous functions was proven by Brouwer and it has been used across numerous fields of mathematics (see \cite{Bo85}).
\vskip.3truecm
\noindent{\bf Brouwer fixed point Theorem.} 
\emph{\ Every continuous function $\varphi$ from a nonempty convex compact subset $C\subseteq\R^n$ to $C$ itself has a fixed point.}
\vskip.3truecm
A selection of a set-valued map $\Phi:X\rightrightarrows Y$ is a function $\varphi:X\rightarrow Y$ that satisfies $\varphi(x)\in\Phi(x)$ for each $x\in X$.
The Axiom of Choice guarantees that set-valued maps with nonempty values always admit selections, but they may have no additional
useful properties.
Michael \cite{Mi56} proved a series of theorems on the existence of continuous selections that assume the condition of lower semicontinuity of set-valued maps.
We present here only one result \cite[Theorem 3.1$^{\prime\prime\prime}$ (b)]{Mi56}.
\vskip.3truecm
\noindent{\bf Michael selection Theorem.}
\emph{\ Every lower semicontinuous set-valued map $\Phi$ from a metric space to $\R^n$ with nonempty convex values admits a continuous selection.}
\vskip.3truecm

\begin{remark}
The Michael selection Theorem holds more in general when the domain of $\Phi$ is a perfectly normal space.
\end{remark}

Collecting the Brouwer fixed point Theorem and the Michael selection Theorem, we deduce the following fixed point result for lower semicontinuous set-valued maps.

\begin{corollary}\label{cor:fixed point}
Every lower semicontinuous set-valued map $\Phi$ from a nonempty convex compact subset $C\subseteq\R^n$ to $C$ itself with nonempty convex values has a fixed point.
\end{corollary}

Notice that, unlike the famous Kakutani fixed point Theorem (see \cite{Bo85}) in which the closedness of $\gph\Phi$ is required, in Corollary \ref{cor:fixed point} the lower semicontinuity of the set-valued map is needed.
No relation exists between the two results as the following example shows.

\begin{example}
The set-valued map $\Phi:[0,3]\rightrightarrows [0,3]$
\begin{displaymath}
\Phi(x):=\left\{\begin{array}{ll}
\{1\} & \mbox{ if } 0\leq x\leq 1\\
(1,2) & \mbox{ if } 1<x<2\\
\{2\} & \mbox{ if } 2\leq x\leq 3
\end{array}\right.
\end{displaymath}
is lower semicontinuous and the nonemptiness of $\fix\Phi$ is guaranteed by Corollary \ref{cor:fixed point}.
Notice that $\fix\Phi=[1,2]$.
Nevertheless the Kakutani fixed point Theorem does not apply since $\gph\Phi$ is not closed.

On the converse, the set-valued map $\Phi:[0,3]\rightrightarrows [0,3]$
\begin{displaymath}
\Phi(x):=\left\{\begin{array}{ll}
\{1\} & \mbox{ if } 0\leq x<1\\
{}[1,2] & \mbox{ if } 1\leq x\leq 2\\
\{2\} & \mbox{ if } 2<x\leq 3
\end{array}\right.
\end{displaymath}
has closed graph and the nonemptiness of $\fix\Phi$ is guaranteed by the Kakutani fixed point Theorem.
Again $\fix\Phi=[1,2]$.
Since $\Phi$ is not lower semicontinuous, Corollary \ref{cor:fixed point} can not be applied.
\end{example}

We conclude this section recalling some topological notations.
Given two subsets $A\subseteq C\subseteq\R^n$ we denote by $\inte_C A$ and $\cl_C A$ the interior and the closure of $A$ in the relative topology of $C$ while $\partial_C A$ indicates the boundary of $A$ in $C$, i.e.
$$\partial_C A:=\cl_C A\setminus \inte_CA = \cl_C A\cap \cl_C (C\setminus A).$$
Lastly $C$ is connected if and only if the subsets of $C$ which are both open and closed in $C$ are $C$ itself and the empty set.

\section{Existence results}

From now on, $C\subseteq \R^n$ is a nonempty convex compact set and $f:C\times C\rightarrow \R$ is an equilibrium bifunction, that is $f(x,x)=0$ for all $x\in C$.
The equilibrium problem is defined as follows:
\begin{equation}\label{eq:ep}
\mbox{find } x\in C \mbox{ such that } f(x,y)\ge 0\mbox{ for all } y\in C.
\end{equation}
Equilibrium problem has been traditionally studied assuming that $f$ is upper semicontinuous in its first argument and quasiconvex in its second one.
Under such assumptions, the issue of sufficient conditions for existence of solutions of (\ref{eq:ep}) was the starting point in the study of the problem.
Ky Fan \cite{Fa72} proved a famous minimax inequality assuming compactness of $C$ and his result holds in a Hausdorff topological vector space.
However, there is the possibility to slightly relax the continuity condition when the vector space is finite dimensional.
The set-valued map
\begin{equation}\label{eq:mapF}
F(x):=\{y\in C:f(x,y)<0\}
\end{equation}
defined on $C$ plays a fundamental role in the formulation of our results.
Clearly $F$ has open lower sections and convex values under the Ky Fan assumptions on the bifunction $f$, that is upper semicontinuity with respect to the first variable and quasiconvexity with respect to the second one.
The fact that $F$ has open lower sections implies that $F$ is lower semicontinuous.
If $F$ had nonempty values, Corollary \ref{cor:fixed point} guarantees the existence of a fixed point of $F$.
This contradicts the fact that $f(x,x)\geq 0$.
Therefore there exists at least one $\bar x$ such that $F(\bar x)=\emptyset$, that is a solution of the equilibrium problem (\ref{eq:ep}).
The following result holds.
\vskip.3truecm
\noindent{\bf Ky Fan minimax inequality.}
\emph{A solution of (\ref{eq:ep}) exists whenever the set-valued map $F$ given in (\ref{eq:mapF}) is lower semicontinuous and convex-valued.}
\vskip.3truecm
After describing this auxiliary result, we focus on the main aim of the paper.
A quasiequilibrium problem is an equilibrium problem in which the constraint set is subject to modifications depending on the considered point.
This format reads
\begin{equation}\label{eq:qep}
\mbox{find } x\in K(x) \mbox{ such that } f(x,y)\ge 0\mbox{ for all } y\in  K(x),
\end{equation}
where $K:C\rightrightarrows C$ is a set-valued map.
Our first existence result is the following.

\begin{theorem}\label{th:existenceQEP}
Assume that $K$ is lower semicontinuous with nonempty convex values and $\fix K$ is closed.
Moreover suppose that
\begin{enumerate}
\item[\em i)]
$F$ is convex-valued on $\fix K$,
\item[\em ii)]
$F$ is lower semicontinuous on $\fix K$,
\item[\em iii)]
$F\cap K$ is lower semicontinuous on $\partial_C\fix K$,
\end{enumerate}
where $F$ is the set-valued map given in (\ref{eq:mapF}).
Then the quasiequilibrium problem (\ref{eq:qep}) has a solution.
\end{theorem}
{\bf Proof.} Corollary \ref{cor:fixed point} ensures the nonemptiness of $\fix K$.
If $\fix K=C$, the existence of solutions to the quasiequilibrium problem descends from the above mentioned Ky Fan minimax inequality.
Otherwise, since $\fix K$ is closed and $\partial_C\fix K =\fix K\setminus \inte_C\fix K$, the emptiness of $\partial_C\fix K$ it would be equivalent to $\fix K=\inte_C\fix K$.
Therefore $\fix K$ would be both open and closed in $C$.
Since every convex set is connected, the only nonempty open and closed subset of $C$ is $C$ itself and this contradicts the fact that $\fix K\ne C$.

Assume that $\inte_C\fix K\ne\emptyset$ (the case $\inte_C\fix K=\emptyset$ is similar and will be shortly discussed at the end of the proof) and define $G:C\rightrightarrows C$ as follows
\begin{displaymath}
G(x):=\left\{\begin{array}{ll}
F(x) & \mbox{ if } x\in \inte_C\fix K\\
F(x)\cap K(x) & \mbox{ if } x\in \partial_C\fix K\\
K(x) & \mbox{ if } x\notin \fix K
\end{array}\right.
\end{displaymath}
The proof is complete if we can show that $G(x)=\emptyset$ for some $x\in C$.
Indeed, since $K$ has nonempty values, then $x\in\fix K$ and two cases are possible.
If $x\in\partial_C\fix K$, then it solves (\ref{eq:qep}); if $x\in\inte_C\fix K$ then it solves (\ref{eq:ep}).
In both cases the quasiequilibrium problem has a solution.

Assume by contradiction that $G$ has nonempty values.
Next step is to prove the lower semicontinuity of $G$.
Fix $x\in C$ and an open set $\Omega\subseteq\R^n$ such that $G(x)\cap\Omega\cap C\ne\emptyset$.
We distinguish three cases.
\begin{enumerate}
\item[a)]
If $x\in\inte_C\fix K$, from the lower semicontinuity of $F$ there exists a neighborhood $U'_x$ such that
\begin{displaymath}
F(x')\cap\Omega\cap C\ne\emptyset,\qquad\forall x'\in U'_x\cap\fix K
\end{displaymath}
which implies
\begin{displaymath}
G(x')\cap\Omega\cap C\ne\emptyset,\qquad\forall x'\in U'_x\cap\inte_C\fix K.
\end{displaymath}
Since $U'_x\cap\inte_C\fix K$ is open in $C$, then $G$ is lower semicontinuous at $x$.
\item[b)]
If $x\in\partial_C\fix K=\partial_C(C\setminus\fix K)$ from the lower semicontinuity of $F$, $K$ and $F\cap K$ there exist neighborhoods $U'_x$, $U''_x$ and $U'''_x$  such that
\begin{eqnarray*}
F(x')\cap\Omega\cap C\ne\emptyset, & \qquad & \forall x'\in U'_x\cap\fix K,\\
K(x')\cap\Omega\cap C\ne\emptyset, & \qquad & \forall x'\in U''_x\cap C,\\
F(x')\cap K(x')\cap\Omega\cap C\ne\emptyset, & \qquad & \forall x'\in U'''_x\cap \partial_C\fix K.
\end{eqnarray*}
Then
\begin{displaymath}
G(x')\cap\Omega\cap C\ne\emptyset,\qquad\forall x'\in U'_x\cap U''_x\cap U'''_x\cap C,
\end{displaymath}
i.e. $G$ is lower semicontinuous at $x$.
\item[c)]
Finally, if $x\notin\fix K$, from the lower semicontinuity of $K$ there exists a neighborhood $U'_x$ such that
\begin{displaymath}
K(x')\cap\Omega\cap C\ne\emptyset,\qquad\forall x'\in U'_x\cap C.
\end{displaymath}
Then
\begin{displaymath}
G(x')\cap\Omega\cap C\ne\emptyset,\qquad\forall x'\in U'_x\cap (C\setminus\fix K).
\end{displaymath}
Since $U'_x\cap(C\setminus\fix K)$ is open in $C$, then $G$ is lower semicontinuous at $x$.
\end{enumerate}
Since by assumption $G$ is also convex-valued, then all the conditions of Corollary \ref{cor:fixed point} are satisfied and there exists $x\in\fix G$.
Clearly $x\in\fix K$ and therefore $x\in\fix F$ which implies $f(x,x)<0$ and contradicts the assumption on $f$.

The issue of $\inte_C\fix K=\emptyset$ remains to be seen.
In this case $\partial_C\fix K=\cl_C\fix K=\fix K$ and $G$ assumes the following form
\begin{displaymath}
G(x):=\left\{\begin{array}{ll}
F(x)\cap K(x) & \mbox{ if } x\in\fix K\\
K(x) & \mbox{ if } x\notin \fix K
\end{array}\right.
\end{displaymath}
The result is obtained by adapting the argument used before.
\qed

\begin{remark}
It is clear from the proof that the assertion remains valid if $f(x,x)=0$ on $C\times C$ is replaced by the weaker $f(x,x)\ge0$ for all $x\in\fix K$.
\end{remark}

\begin{remark}\label{re:alternative}
The proof of Theorem \ref{th:existenceQEP} allows to establish that a solution of (\ref{eq:qep}) belongs to
\begin{displaymath}
\partial_C\fix K\cup \{x\in \inte_C\fix K: x \mbox{ solves } (\ref{eq:ep})\}.
\end{displaymath}
In particular if (\ref{eq:ep}) has no solution then Theorem \ref{th:existenceQEP} ensures that a solution of (\ref{eq:qep}) lies on the boundary of $\fix K$.
\end{remark}

\begin{remark}\label{re:fan}
By specializing to $K(x):=C$, for all $x\in C$, Theorem \ref{th:existenceQEP} becomes the Ky Fan minimax inequality.
Indeed $\fix K=C$ and conditions i) and ii) coincide with the assumptions in Ky Fan minimax inequality.
Instead, since $\partial_C\fix K=\emptyset$, condition iii) is trivially satisfied.
\end{remark}

Theorem \ref{th:existenceQEP} is strongly related to \cite[Lemma 3.1]{Cu95}.
The two sets of conditions differ only in that the lower semicontinuity of $F\cap K$ on the whole space $C$ assumed in \cite[Lemma 3.1]{Cu95} is here replaced by the lower semicontinuity of $F$ on $\fix K$ and the lower semicontinuity of $F\cap K$ on $\partial_C \fix K$.
We provide an example in which the results are not comparable to each other.

\begin{example}
Let $C:=[0,1]$ and
\begin{displaymath}
f(x,y):=\left\{\begin{array}{ll}
-1 & \mbox{ if } x=0 \mbox{ and } y\in(0,1]\\
0 & \mbox{ otherwise}
\end{array}\right.
\end{displaymath}
If $K(x):=\{x\}$, for all $x\in [0,1]$, then $F\cap K=\emptyset$ is trivially lower semicontinuous and the assumptions of \cite[Lemma 3.1]{Cu95} are satisfied.
Instead $F$ is not lower semicontinuous at $0\in\fix K=[0,1]$.

On the other hand if $K(x):=\{1-x\}$, for all $x\in [0,1]$, then $\fix K=\{1/2\}$, the assumptions of Theorem \ref{th:existenceQEP} are trivially satisfied,
but $F\cap K$ is not lower semicontinuous at $0$.
\end{example}

It would be desirable to find more tractable conditions on $f$, disjoint from the ones assumed on $K$, which guarantee that all the assumptions  i), ii) and iii) of Theorem \ref{th:existenceQEP} are satisfied.
Clearly the convexity of $F(x)$ can be deduced from the quasiconvexity of $f(x,\cdot)$ for all $x\in\fix K$.
While the upper semicontinuity of $f(\cdot,y)$ on $\fix K$ implies that $F^{-1}(y)$ is open on $\fix K$ and hence $F$ is lower semicontinuous on $\fix K$.

The last part of this section is devoted to furnish sufficient conditions for assumption iii), i.e. which guarantee the lower semicontinuity of the set-valued map $F\cap K$ on $\partial_C\fix K$.
We propose two approaches.
The former one consists in exploiting the following result in \cite{Pa91}.

\begin{proposition}\label{pr:lsc intersection}
Let $\Phi_1,\Phi_2:X\rightrightarrows Y$ be set-valued maps between two topological spaces.
Assume that $\gph\Phi_1$ is open on $X\times Y$ and $\Phi_2$ is lower semicontinuous.
Then $\Phi_1\cap\Phi_2$ is lower semicontinuous.
\end{proposition}

Since $K$ is assumed to be lower semicontinuous, we investigate which assumptions ensure the open graph of $F$ given in (\ref{eq:mapF}),
that is the openness of the set
\begin{equation}\label{eq:open}
\{(x,y)\in \partial_C\fix K\times C:f(x,y)< 0\}.
\end{equation}
Hence, Theorem \ref{th:existenceQEP} still works by using this condition instead of iii).
It is interesting to compare this fact with \cite[Theorem 2.1]{Cu97} where the openness of the set $\{(x,y)\in C\times C:f(x,y)< 0\}$ is required instead of the openness of (\ref{eq:open}) and the lower semicontinuity of $F$ on $\fix K$.
One should not overlook the fact that even though the results are formally similarly formulated, unlike our result, \cite[Theorem 2.1]{Cu97} does not reduce to Ky Fan minimax inequality when $K(x)=C$, for all $x\in C$.

An open graph result is \cite[Proposition 2]{Zh95} which affirms that if $X$ is a topological space and $\Phi:X\rightrightarrows \R^n$ is a set-valued map with convex values, then $\Phi$ has open graph in $X\times\R^n$ if and only if $\Phi$ is lower semicontinuous and open valued.
This fact has been used to establish the existence of continuous selections, maximal elements, and fixed points of correspondences in various economic applications.

Up to translations, this result also holds when the codomain of $\Phi$ is an affine subset of $\R^n$ \cite[Theorem 1.12]{Yu98}.
We recall that an affine set of $\R^n$ is the translation of a vector subspace.
Moreover, the affine hull of a set $C$ in $\R^n$, which is denoted by $\aff C$, is the smallest affine set containing $C$, or equivalently, the intersection of all affine sets containing $C$.

\begin{theorem}\label{th:sufficientconditions1}
Let $A\supseteq C$ be an open set on $\aff C$ and $\hat{f}:C\times A\rightarrow \R$ be a bifunction such that $\hat{f}(x,y)=f(x,y)$ for all $(x,y)\in C\times C$.
Denote by $\hat{F}$ the set-valued map
\begin{displaymath}
\hat{F}(x):=\{y\in A:\hat{f}(x,y)<0\}
\end{displaymath}
defined on $C$ and assume that $K$ is lower semicontinuous with nonempty convex values and $\fix K$ is closed.
Moreover suppose that
\begin{enumerate}
\item[\em i)]
$\hat{F}$ is convex-valued on $\fix K$,
\item[\em ii)]
$\hat{F}$ has open lower sections on $\fix K$,
\item[\em iii)]
$\hat{F}(x)$ is open on $\aff C$ for all $x\in\partial_C\fix K$.
\end{enumerate}
Then the quasiequilibrium problem (\ref{eq:qep}) has a solution.
\end{theorem}
{\bf Proof.} We have to show that all the assumptions of Theorem \ref{th:existenceQEP} are fulfilled.
Since the set-valued map $F$ given in (\ref{eq:mapF}) can be expressed as $\hat F\cap C$, i) implies that $F$ is convex-valued on $\fix K$ and ii) implies that $F$ is open lower section on $\fix K$.
In particular $F$ is lower semicontinuos on $\fix K$.
Furthermore assumption iii) allows to apply \cite[Theorem 1.12]{Yu98} which ensures that $\gph \hat F$ is open on $\partial_C\fix K\times \aff C$.
Hence $\gph F=\gph \hat F\cap (\partial_C\fix K\times C)$ is open on $\partial_C\fix K\times C$ and
Proposition \ref{pr:lsc intersection} guarantees that the intersection map $F\cap K$ is lower semicontinuous on $\partial_C\fix K$.
\qed

The open graph result \cite[Proposition 2]{Zh95} no longer holds when $\R^n$ (or an affine space) is replaced with an infinite dimensional Hilbert space \cite{Ba12}.
However if $C\subset\R^n$ is a polytope, that is the convex hull of a finite set, then every $\Phi:X\rightrightarrows C$ with open lower sections and convex open values has open graph \cite[Proposition 11.14]{Bo85}.
This fact can be used for proving our next result.

\begin{theorem}\label{th:sufficientconditions2}
Assume that $C$ is a polytope and $K$ is lower semicontinuous with nonempty convex values and $\fix K$ is closed.
Moreover suppose that
\begin{enumerate}
\item[\em i)]
$F$ is convex-valued on $\fix K$,
\item[\em ii)]
$F$ has open lower sections on $\fix K$,
\item[\em iii)]
ì$F(x)$ is open on $C$ for all $x\in \partial_C \fix K$,
\end{enumerate}
where $F$ is the set-valued map given in (\ref{eq:mapF}).
Then the quasiequilibrium problem (\ref{eq:qep}) has a solution.
\end{theorem}
{\bf Proof.} The set-valued map $F$ has open lower sections, convex and open values.
Then its graph is open on $\partial_C \fix K\times C$ \cite[Proposition 11.14]{Bo85} and the lower semicontinuity of $F\cap K$ follows from Proposition \ref{pr:lsc intersection}.
\qed

Notice that the lower semicontinuity condition ii) assumed in Theorem \ref{th:existenceQEP} has been replaced in the last two results by the requirement that the lower sections are open.
This is due to two different reasons.

In the proof of Theorem \ref{th:sufficientconditions1}, in order to apply \cite[Theorem 1.12]{Yu98} and get that $\gph\hat F$ is open, it would be enough to require the lower semicontinuity of $\hat F$.
However such an assumption would not guarantee the lower semicontinuity of $F=\hat F\cap C$ which is assumption ii) in Theorem \ref{th:existenceQEP}.

On the other hand, assumption ii) in Theorem \ref{th:sufficientconditions2} is necessary to get the openness of $\gph F$ as a consequence of \cite[Proposition 11.14]{Bo85}.
The next example shows that a set-valued map $\Phi$ acting from a topological vector space to a polytope $C$ may not have open graph and \cite[Proposition 11.14]{Bo85} fails even if it is lower semicontinuous with convex and open values.

\begin{example}
Let $C:=\{(x,y)\in \R^2:|x|+|y|\leq 1\}$ be a closed convex set in $\R^2$.
The set-valued map $\Phi:[0,1]\rightarrow C$ defined by
\begin{displaymath}
\Phi(t):=\left\{\begin{array}{ll}
C\setminus \{(x,y):x+y=1\} & \mbox{ if } t>0\\
C & \mbox{ if } t=0
\end{array}\right.
\end{displaymath}
is lower semicontinuous with convex open values in $C$ but it has not open lower sections since $\phi^{-1}(0,1)=\{0\}$.
Nevertheless $\gph\Phi$ is not open in $[0,1]\times C$ since the sequence $ \{(n^{-1},1-n^{-1},n^{-1})\}\in [0,1]\times C$ does not belong to $\gph\Phi$ but its limit $(0,1,0)\in\gph\Phi$.
\end{example}

We answer in the negative the question posed in \cite{BePaRa76} where the authors affirm that they do not know whether \cite[Proposition 11.14]{Bo85} can be generalized to the case where $C$ is an arbitrary convex subset of $\R^n$.
This also explains why we need to extend the domain of $f(x,\cdot)$ from $C$ to an open subset of $\aff C$ in Theorem \ref{th:sufficientconditions1}.
\begin{example}\label{ex:graphnotopen}
Let $C\subseteq \R^2$ be the closed unit ball.
The set-valued map $\Phi:[0,1]\rightrightarrows C$ defined by
\begin{displaymath}
\Phi(x):=\left\{\begin{array}{ll}
C\setminus \{(\cos x,\sin x)\} & \mbox{ if } x>0\\
C & \mbox{ if } x=0
\end{array}\right.
\end{displaymath}
has open lower sections and convex open values in $C$.
Nevertheless $\gph\Phi$ is not open in $[0,1]\times C$.
Indeed $(1,0)\in\Phi(0)$ and there is no neighborhood $U$ of $(1,0)$ such that $U\cap C\subseteq\Phi(x)$ for $x$ small enough.
\end{example}

A second possible approach for the lower semicontinuity of $F\cap K$ could be to show the nonemptiness of the intersection between the interior of $F$ and $K$.
Indeed \cite[Corollary 1.3.10]{BoGeMyOb84} affirms that the set-valued map $\Phi_1\cap\Phi_2$ is lower semicontinuous on the topological space $X$ provided that $\Phi_1,\Phi_2:X\rightrightarrows C$ are convex-valued, lower semicontinuous set-valued maps and
\begin{equation}\label{eq:inte}
\Phi_1(x)\cap\Phi_2(x)\neq \emptyset\quad\Rightarrow\quad\Phi_1(x)\cap\inte\Phi_2(x)\neq \emptyset.
\end{equation}
The following example shows that such result could not be guaranteed (as erroneously stated in \cite[Theorem 1.13]{Yu98}) if the interior is replaced by the relative interior in condition (\ref{eq:inte}).
Given a set $C\subseteq\R^n$, we denote by $\ri C$ the relative interior of $C$, namely, $\ri C=\inte_{\aff C}C$.

\begin{example}
Let $C\subseteq \R^2$ be the closed unit ball and $\Phi_1:[0,1]\rightrightarrows C$ be defined as in Example \ref{ex:graphnotopen}.
Consider $\Phi_2:[0,1]\rightrightarrows C$ defined by
\begin{displaymath}
\Phi_2(x):=\{(\cos x,\sin x)\}\qquad \forall x\in[0,1].
\end{displaymath}
Then $\Phi_2$ is a continuous single-valued map and $\Phi_1$ is convex-valued with open lower sections.
Furthermore
\begin{displaymath}
\Phi_1(x)\cap\Phi_2(x)=\left\{\begin{array}{ll}
\emptyset & \mbox{ if } x>0\\
\{(1,0)\} & \mbox{ if } x=0
\end{array}\right.
\end{displaymath}
and $\Phi_1(0)\cap \ri\Phi_2(0)=C\cap \{(1,0)\}=\{(1,0)\}$.
Nevertheless $\Phi_1\cap\Phi_2$ is not lower semicontinuous at $0$.
Notice that $\Phi_1(x)$ is even open on $C$, for all $x\in[0,1]$.
\end{example}

The following is a correct version of \cite[Theorem 1.13]{Yu98}.

\begin{proposition}\label{pr:lsc intersection1}
Let $X$ be a topological space, $C\subseteq\R^n$ and $\Phi_1,\Phi_2:X\rightrightarrows C$ be lower semicontinuous and convex-valued.
Moreover, for all $x\in X$ assume that $\aff \Phi_2(x)=\aff C$ and
\begin{displaymath}
\Phi_1(x)\cap\Phi_2(x)\neq \emptyset\quad \Rightarrow\quad\Phi_1(x)\cap\ri\Phi_2(x)\neq \emptyset
\end{displaymath}
then $\Phi_1\cap \Phi_2$ is lower semicontinuous.
\end{proposition}
{\bf Proof.} By definition, up to isomorphism, there exists $m\leq n$ such that $\aff C=x_0+\R^m$, where $x_0\in C$ is arbitrarily fixed.
Define $\hat\Phi_i:X\rightrightarrows \R^m$ by
$\hat\Phi_i:=\Phi_i-x_0$, $i=1,2$.
Then $\hat\Phi_1$ and $\hat\Phi_2$ are lower semicontinuous and convex-valued.
Furthermore, since $\aff\Phi_2(x)=\aff C$, then $\ri\Phi_2(x)=x_0+\inte\hat\Phi_2(x)$ and
$\hat\Phi_1(x)\cap \inte \hat\Phi_2(x)\neq \emptyset$ whenever $\hat\Phi_1(x)\cap \hat\Phi_2(x)\neq \emptyset$.
By \cite[Corollary 1.3.10]{BoGeMyOb84} it follows that $\hat\Phi_1\cap \hat\Phi_2$ is lower semicontinuous.
This means in turn that $\Phi_1\cap \Phi_2$ is lower semicontinuous.
\qed

Now we are in position to prove our last existence result.

\begin{theorem}\label{th:sufficientconditions3}
Assume that $K$ is lower semicontinuous with nonempty convex values and $\fix K$ is closed.
Moreover suppose that
\begin{enumerate}
\item[\em i)]
$F$ is convex-valued on $\fix K$,
\item[\em ii)]
$F$ is lower semicontinuous on $\fix K$,
\item[\em iii)]
$\aff K(x)=\aff C$, for all $x\in \partial_C \fix K$,
\item[\em iv)]
$F(x)$ is open on $C$, for all $x\in \partial_C \fix K$,
\end{enumerate}
where $F$ is the set-valued map given in (\ref{eq:mapF}).
Then the quasiequilibrium problem (\ref{eq:qep}) has a solution.
\end{theorem}
{\bf Proof.} It is enough to show that assumption iii) of Theorem \ref{th:existenceQEP} holds, i.e.
$F\cap K$ is lower semicontinuous on $\partial_C\fix K$.
Let $x\in\partial_C\fix K$ be fixed and assume that $F(x)\cap K(x)\neq\emptyset$
(otherwise the intersection is trivially lower semicontinuous at $x$).
By assumption there exists an open set $\Omega\subseteq \R^n$ such that $F(x)=\Omega\cap C$.
Then
\begin{displaymath}
\emptyset\neq F(x)\cap K(x)=\Omega\cap C\cap K(x)=\Omega\cap K(x).
\end{displaymath}
From \cite[Corollary 6.3.2]{Ro70} we get
\begin{displaymath}
\emptyset\neq \Omega\cap \ri K(x)=F(x)\cap \ri K(x)
\end{displaymath}
The lower semicontinuity of $F\cap K$ at $x$ follows from Proposition \ref{pr:lsc intersection1}.
\qed

Now we make a comparison with an analogous result in \cite{Cu95}.
The assumptions of Theorem \ref{th:sufficientconditions3} are the same as those of \cite[Theorem 3.2]{Cu95} except that conditions iii) and iv) must be verified for all $x\in \partial_C\fix K$ instead of for all $x\in C$.
Thus, Theorem \ref{th:sufficientconditions3} is clearly more general and, unlike \cite[Theorem 3.2]{Cu95}, it reduces to Ky Fan minimax inequality when the constraint set-valued map $K$ is equal to $C$.

\section{Conclusions}

In this paper existence results for the solution of finite dimensional quasiequilibrium problems are obtained by using a Michael selection result for lower semicontinuous set-valued maps.
The peculiarity of our results, which make them different from other results in the literature to the best of knowledge of the authors, is the fact that they reduce to Ky Fan minimax inequality when the constraint map is constant.

Moreover we provide information regarding the position of a solution.
In fact either it is a fixed point of the constraint set-valued map which solves an equilibrium problem or it lies in the boundary of the fixed points set.
To know this property seems promising for the construction of solution methods.
Future works could be devoted to exploit such result to propose computational techniques for solving quasiequilibrium problems.

Another possible advance consists in studying conditions which permit to replace the compactness of the domain with suitable coercivity conditions on the equilibrium bifunction.


\begin{thebibliography}{00}

\bibitem{Fa72}
Fan K.:
A minimax inequality and applications.
In: Shisha O. (ed.): Inequalities III, pp. 103--113. Academic Press, New York (1972)

\bibitem{BiCaPaPa13}
Bigi G., Castellani M., Pappalardo M., Passacantando M.:
Existence and solution methods for equilibria.
European J. Oper. Res. 227, 1--11 (2013)

\bibitem{BeGoLi73}
Bensoussan A., Goursat M., Lions J.L.:
Contr\^{o}le impulsionnel et in\'{e}quations quasi-variationnelles stationnaires.
C.R. Acad. Sci. Paris S\'{e}r. A 276, 1279--1284 (1973)

\bibitem{Mo76}
Mosco U.:
Implicit variational problems and quasi variational inequalities.
In: Lecture Notes in Math., vol. 543, pp. 83--156. Springer-Verlag, Berlin (1976)

\bibitem{AlRa16}
Alleche B., R\u{a}dulescu, V.D.:
Solutions and approximate solutions of quasi-equilibrium problems in Banach spaces.
J. Optim. Theory Appl. 170, 629--649 (2016)

\bibitem{Au93}
Aubin J.P.:
Optima and equilibria.
Springer-Verlag, Berlin (1993)

\bibitem{AuCoIu17}
Aussel D., Cotrina J., Iusem A.:
Existence results for quasi-equilibrium problems.
J. Convex Anal. 24, 55--66 (2017)

\bibitem{CaGi15}
Castellani M., Giuli M.:
An existence result for quasiequilibrium problems in separable Banach spaces.
J. Math. Anal. Appl. 425, 85--95 (2015)

\bibitem{CaGi16}
Castellani M., Giuli M.:
Approximate solutions of quasiequilibrium problems in Banach spaces.
J. Global Optim. 64, 615--620 (2016)

\bibitem{Cu95}
Cubiotti P.:
Existence of solutions for lower semicontinuous quasiequilibrium problems.
Comput. Math. Appl. 30, 11--22 (1995)

\bibitem{Cu97}
Cubiotti P.:
Existence of Nash equilibria for generalized games without upper semicontinuity.
Internat. J. Game Theory 26, 267--273 (1997)

\bibitem{Mi56}
Michael E.:
Continuous selections. I.
Ann. of Math. 63, 361--382 (1956)

\bibitem{Bo85}
Border K.C.:
Fixed point theorems with applications to economics and game theory.
Cambridge University Press, Cambridge (1985)

\bibitem{Pa91}
Papageorgiou N.S.:
On the existence of $\psi$-minimal viable solutions for a class of differential inclusions.
Arch. Math. 27, 175--182 (1991)

\bibitem{Zh95}
Zhou J.:
On the existence of equilibrium for abstract economies,
J. Math. Anal. Appl. 193, 839--858 (1995)

\bibitem{Yu98}
Yuan G.X.-Z.:
The study of minimax inequalities and applications to economies and variational inequalities.
Memoirs of the American Mathematical Society, vol. 132.
Providence, Rhode Island (1998)

\bibitem{Ba12}
Bagh A.:
Lower hemi-continuity, open sections, and convexity: counter examples in infinite dimensional spaces.
Theoret. Econom. Lett. 2, 121--124 (2012)

\bibitem{BePaRa76}
Bergstrom T.C., Parks R.P., Rader T.:
Preferences which have open graphs.
J. Math. Econom. 3, 265--268 (1976)

\bibitem{BoGeMyOb84}
Borisovich Y., Gel'man B.D., Myshkis A.D., Obukhovskii V.V.:
Multivalued  mappings.
J. Soviet Math. 24, 719--791 (1984)

\bibitem{Ro70}
Rockafellar R.T.:
Convex Analysis.
Princeton University Press, Princeton (1970)
\end{thebibliography}
\end{document}